\magnification=\magstep1
\hfuzz=2pt
\headline={\ifnum\pageno=1 \hfil\else\hss{\tenrm\folio}\hss\fi}
\footline={\hfil}
\font\titlefont= cmbx10 scaled \magstep3

\def\integer{{\bf N}}
\def\real{{\bf R}}

\def\Romannumeral#1{\uppercase\expandafter{\romannumeral#1}}
\def\date{\line{\number\day/\number\month/\number\year\hfil}}

\def\expectation{\mathchoice{\rm I\hskip-1.9pt E}{\rm I\hskip-1.9pt E}
{\rm I\hskip-.8pt E}{\rm I\hskip-1.9pt E}}

\def\Oun{{\cal O}(1)}

\def\proof{\noindent{\bf Proof. }}

\font\bigmath=cmmi10 scaled \magstep2
\def\bigchi{\hbox{\bigmath \char31}}

\newskip\refskip\refskip=4em
\def\refsize{\advance\leftskip by \refskip}
\def\ref#1#2{\noindent\hskip -\refskip\hbox to
\refskip{[#1]\hfil}{\noindent #2\hfil}\medskip}
\def\proba{\mathchoice{\rm I\hskip-1.9pt P}{\rm I\hskip-1.9pt P} 
{\rm I\hskip-.9pt P}{\rm I\hskip-1.9pt P}}

\vglue 3cm
\centerline{\titlefont STATISTICS OF CLOSEST RETURN}
\vskip 1cm 
\centerline{\titlefont FOR SOME NON UNIFORMLY}
\vskip 1cm
\centerline{\titlefont HYPERBOLIC SYSTEMS.}
\vskip 3cm
\centerline{P.COLLET}
\bigskip
\centerline{Centre de Physique Th\'eorique}
\centerline{Laboratoire CNRS UMR 7644}
\centerline{Ecole Polytechnique}
\centerline{F-91128 Palaiseau Cedex (France)}
\vskip 3cm
\noindent{\bf Abstract}: 
For non uniformly hyperbolic maps  of the interval with exponential
decay of correlations we prove that 
the law of closest return to a given point when suitably normalized
is almost surely
asymptotically exponential. A similar result holds when the reference
point is the initial point of the trajectory. We use the framework for
non uniformly hyperbolic dynamical systems 
developed by L.S.Young.   
\vfill\supereject

\beginsection{I. INTRODUCTION.}

The statistics of entrance time in a small set is a long standing
important problem. One of the first instances is the famous question
raised by Boltzmann about the time it takes to observe the clustering of
all the molecules of a gas in only half of the available volume. Other
important applications include the occurrence of rare events, eventually
catastrophic.  More recently, several algorithms were proposed to
measure various quantities like dimension and entropy of dynamical
systems which involve such entrance time questions. We refer to [ABST]
for a review of these algorithms. Another application concerns the
optimal compression of data sets. The well known compression algorithm
developed by Ziv and Wyner is based on the coding of repetitions of
patterns previously coded. Its optimality is based on almost sure
results for the typical recurrence time which were proven for general
ergodic sources by Ornstein and Weiss. We refer to [WZW] and [Sh] for
recent reviews on this subject and references to older works. This
algorithm can also be viewed as a way of measuring the entropy.

In all these questions, the asymptotic result is determined by a law of
large numbers. One would like to understand the fluctuations in order to
control the rate of convergence and for statistical purposes.  The rate
of convergence for the Wyner-Ziv algorithm was recently derived for the
case of sufficiently mixing sources in [CGS], [K] and [P].

In the present paper we derive similar results from a  different
point of view which emphasizes the topology of the phase space, 
and for systems which are non uniformly hyperbolic. The problem
discussed below can be formulated in terms of extreme statistics (see
[G]). Consider a dynamical system given by a (compact) phase space $\Omega$
equipped with a metric $d$, a
continuous map $T$ on $\Omega$ and a $T$ invariant ergodic probability 
measure $\mu$.
 Assume a point $x$ in phase space has been chosen 
and define a sequence $(X_{j})$ of real random variables on $\Omega$ by 
$$
X_{j}(y)=-\log(d(x,T^{j}(y)))\;.
$$ 
Let $(Z_{n})$ be the sequence of successive maxima of the sequence
$(X_{j})$, namely
$$
Z_{n}(y)=\sup_{0\le j\le n}X_{j}(y)\;.
$$
One may conjecture that if the dynamical system satisfies the
Eckmann-Ruelle conjecture, then almost surely
$$
\lim_{n\to\infty}{Z_n\over \log n}=D^{-1}
$$
where $D$ is the dimension of the measure. A more precise question is to
  to ask if there are two sequences of positive
numbers $(a_{n})$ and $(b_{n})$ such that for any positive real number
$s$, the sequence 
$$
\proba\left(Z_{n}>a_{n}s+b_{n}\right)
$$
converges. We will denote by $\proba$ the probability associated to the
measure $\mu$ and by $\expectation$ the corresponding expectation.
We refer to [G] for the study of the statistics of extremes for
 independent random
variables and for related questions. In the present case of dynamical
systems, one may expect that the result (and in particular the choice of
the two sequences $(a_{n})$ and $(b_{n})$)  depends on the point $x$
chosen at the beginning. We will assume below that the point $x$ has
been chosen at random with respect to $\mu$ and we will prove results
almost surely with respect to this choice. In order to simplify the
notation we will usually not mention the $x$ dependence since this is a
point which is chosen once for all. When this dependence needs to be
emphasized we will denote it by an exponent in $(X_{j}^{x})$ and
$(Z_{n}^{x})$.

We now describe the dynamical systems  for which the result will be proven.
An abstract frame for non uniformly hyperbolic systems 
was introduced by L.S.Young in [Y1] and [Y2] (see also references
therein and [BV] and [KN] for earlier constructions).
Instead of using the completely abstract formulation we will rather keep
the equivalent version in the phase space. The reason for doing so is
that the topology is somewhat obscured when the system is lifted to the
abstract context. We will also work explicitly the case of maps of the
interval with exponential decay of correlations although several results
extend to  more general situations. We will mention some of these
extensions when appropriate. We now formulate the hypothesis on our
dynamical system which follow directly from the work of L.S.Young.

We will consider a $C^{2}$ map $f$ of the interval $[a,b]$ into itself and
we denote by $K$ the sup norm of its derivative
$$
K=\sup_{x\in[a,b]}|f'(x)|\;.
$$
We assume there is an open interval $\Lambda$ in $[a,b]$ with dense
orbit and  with the following properties. 

\item{{\bf H1}} {\sl 
There exists 
 a sequence $(R_{i})_{i\in\integer}$ of positive integers, with largest
common divisor equal to one  and a  sequence of disjoint open  subintervals
$(\Lambda_{i})_{i\in\integer}$ of $\Lambda$ satisfying 
$\lambda(\Lambda\backslash \cup_{i}\overline\Lambda_{i})=0$
such that the following holds. 
For any $j\in\integer$, $f^{R_{j}}$ is a bijection from
$\Lambda_{j}$ to $\Lambda$. 
There exits also an integer valued function $s$
defined on $\Lambda\times\Lambda$
such that for $x$ and $y$ in $\Lambda$, the orbits of $x$ and $y$ follow
 each other up to time $s(x,y)$ in the sense that the corresponding
 orbits under $f^R$ fall in the same $\Lambda_i$.
We assume  $s$ is finite $\lambda\times\lambda$
almost surely, where $\lambda$ is the Lebesgue measure.}

\noindent 
Hypothesis {\bf H1} is of course of Markov type. Note however that the
number $R_{j}$ may not be the first return to $\Lambda$ of
$\Lambda_{j}$. It is a return where the properties of uniform backward
contraction {\bf H2} and uniform distortions {\bf H3} are satisfied. 
We will speak of these returns as "official" returns.

\item{{\bf H2}} {\sl
There is a constant $C>0$ and a number $0<\beta<1$ such
that for any $x$ and $y$ in $\Lambda$, and any $0\le n\le s(x,y)$ we have
$$
|f^{n}(x)-f^{n}(y)|\le C\beta^{s(x,y)-n}\;.
$$}

\item{{\bf H3}} {\sl For any $x$, $y$ in $\Lambda$
 and any $0\le k\le n\le s(x,y)$ we have
$$
\log\left(\prod_{i=k}^{n}{|f'(f^{i}(x))|\over |f'(f^{i}(y))|}\right)
\le C\beta^{s(x,y)-n}\;.
$$}

In [Y1] and [Y2] examples of dynamical systems where given where
these hypothesis are satisfied. In particular unimodal 
maps of the interval with sufficient instability of the critical orbit
satisfy these hypothesis. There are also examples with neutral fixed
points and higher dimensional cases. We will make some comments below
about  these cases. 

After having described the setting in phase space, we now come to the
invariant measure. Let $\lambda$ denote the Lebesgue measure on $[a,b]$.
 The next
hypothesis concerns the random variable $R$ defined on $\Lambda$ by
$R(x)=R_{j}$ if $x\in \Lambda_{j}$. 

\item{{\bf H4}} {\sl The random variable $R$ is integrable with respect to
$\lambda$ (restricted to $\Lambda$)}. 

One of the first results of L.S.Young is that under the above hypotheses,
there is a measure $\mu_{0}$ on $\Lambda$ which is equivalent to
$\lambda$ (more precisely with a density bounded above and bounded below
away from zero) and which is invariant and ergodic for $f^{R}$. This
leads to an invariant ergodic measure $\mu$ for the map $f$ which is
given by
$$
\mu(A)=Z^{-1}\sum_{l}\sum_{j=0}^{R_{l}-1}\mu_{0}\left(\Lambda_{l}\cap f^{-j}
\big(A\cap f^{j}(\Lambda_{l})\big) \right)\;,\eqno(I.1)
$$
with
$$
Z=\sum_{l}R_{l}\mu_{0}\left(\Lambda_{l}\right)=\int R\;d\mu_{0}\;.
$$

One of the main result in [Y1] and [Y2] is a bound on the decay of
correlations for for H\"older continuous. Namely if $g_1$ is H\"older
continuous and $g_2\in L^{\infty}(d\lambda)$, the decay of correlations
$\alpha(\,\cdot\,)$ defined by
$$
\alpha(n)=\left|\int g_1\;g_2\circ f^n d\mu \;-\int g_1\;d\mu\int
  g_2\;d\mu \right|
$$
is related to the behaviour for large $n$ of $\lambda(R>n)$. If this
sequence decays exponentially fast, the same is true for $\alpha$ with
the same decay rate. A similar result holds in the case of polynomial
decay.  In the case of exponential decay of $\lambda(R>n)$, a stronger
version of the decay of correlations was proven in [Y1] which is
analogous to the case of piecewise expanding maps of the interval.
Although this stronger result would slightly simplify some arguments
below, we will not use it.
We can now formulate our main result.

\proclaim{Theorem I.1}. {Assume the hypotheses {\bf H1}-{\bf H4}. Assume
$\lambda(R>k)$ decays exponentially fast with $k$. Then for $\mu$ (or
Lebesgue) almost every $x$ we have
$$
\lim_{n\to\infty}\proba\left(Z_{n}^{x}<s+\log n\right)=e^{-2h(x)e^{-s}}
$$
where $h$ is the density of the absolutely continuous invariant measure
$\mu$.}

This is sometimes called Gumbel's law.  There is an obvious relation
with the entrance time in a ball of radius $e^{-s}/n$ centered at $x$,
i.e. if $Z_{n}<s+\log n$, the orbit has not entered the ball up to time
$n$.

In the next section we will prove some preparatory results, and in
particular we will control the measure of the set of points which recur
too fast. The proof of Theorem I.1 will then be given in section 3,
inspired by the techniques used for extreme statistics.
In section 4, we will establish a fluctuation result for the case where
the initial point is the reference point. 

We mention also that some intermediate 
results proven below where already derived in
explicit situations in
order to prove hypotheses {\bf H1}-{\bf H4} or the decay of 
$\lambda(R>k)$. One of the goal of this paper is to show that the
previously mentioned hypotheses are sufficient to prove the result
without reference to particular constructions. 

\beginsection{II. ESTIMATES FOR THE SET OF RAPIDLY RECURRING POINTS.}

In this section we will estimate the measure of some sets of points with
exceptional behaviour. For later references we start with the following
easy lemma. 

\proclaim{Lemma II.1}. {
$$
\sum_{l\,,\,R_{l}>k}R_{l}\;\lambda(\Lambda_{l})
\le\cases{ 2\sum_{s=k/2}^{\infty}\;\lambda(\{R>s\})\;,&\cr
\sum_{s=k}^{\infty}\;\lambda(\{R>s\})+k\,\lambda(\{R>k\})\;.
&\cr}
$$}

\proof We have indeed for any $q>0$
$$
\sum_{s=q}^{\infty}\lambda(\{R>s\})=\sum_{s=q}^{\infty}
\sum_{l=1}^{\infty}\lambda(\{R>s\}\cap \Lambda_{l})
$$ 
$$
=\sum_{l,\,R_{l}>q}^{\infty}(R_{l}-q)\lambda(\Lambda_{l})\;.
$$
For the first estimate, we take $q=k/2$ and restrict the last sum in the
above equality to the range $R_{l}>k$ which implies $R_{l}-q\ge R_l/2$ and
the result follows. For the second estimate we simply take $q=k$ and
rearrange the equality.

We will need later  an estimate of the $\mu$ measure of sets
with small Lebesgue measure.

\proclaim{Lemma II.2}. {Assume exponential decay in $k$ of $\lambda(\{R>k\})$.
Then there are two positive constants $C$ and $\theta$ such that
for any Lebesgue measurable set $I$,  we have
$$
\mu(I)\le C\lambda(I)^{\theta} \;.
$$} 

\proof  From formula (I.1), for $0\le j<R_{l}$ we have to consider the sets
$$
I_{l,j}
=I\cap f^{j}(\Lambda_{l})\;.
$$
Since
$f^{j}$ is injective  on $\Lambda_{l}$, there is a set $\tilde I_{j,l}$
in $\Lambda_{l}$ such that $f^{j}(\tilde I_{j,l})=I\cap
f^{j}(\Lambda_{l})$.  There are now two cases. 

In the first case, we assume $K^{R_{l}-j}\le |I|^{-1/2}$. We now use the
distortion bound on $\Lambda_{l}$ for $f^{R_{l}}$. We get
$$
{|\tilde I_{j,l}|\over |\Lambda_{l}|}\le\Oun
{|f^{R_{l}}(\tilde I_{j,l})|\over |f^{R_{l}}(\Lambda_{l})|}
=\Oun {\big|f^{R_{l}-j}\big(I\cap f^{j}(\Lambda_{l})\big)\big|
\over |f^{R_{l}}(\Lambda_{l})|}
\le \Oun K^{R_{l}-j}|I|\le \Oun |I|^{1/2}. 
$$
This implies
$$
|\tilde I_{j,l}|\le \Oun |I|^{1/2} |\Lambda_{l}|\;,
$$
and since $\mu_{0}$ is equivalent to the Lebesgue measure
$$
\mu_{0}(\tilde I_{j,l})\le \Oun |I|^{1/2} \mu_{0}(\Lambda_{l})\;.
$$
We can now sum over $j$ and $l$ to get
$$
\sum_{j,\,l,\, K^{R_{l}-j}< |I|^{-1/2} }\mu_{0}(\tilde I_{j,l})\le
 \Oun |I|^{1/2} \sum_{l}R_{l}\;\mu_{0}(\Lambda_{l})\;.
$$
This last sum is finite since $R$ is integrable with respect to
$\lambda$ by hypothesis {\bf H4}, and $\mu_0$ is equivalent to
$\lambda$.  

We now deal with the second case, namely  $K^{R_{l}-j}> |I|^{-1/2}$.  
This implies of course 
$$
R_{l}\ge {\log |I|^{-1}\over 2\log K}\;. 
$$
Therefore
$$
\sum_{j,\,l,\, K^{R_{l}-j}\ge |I|^{-1/2} }\mu_{0}(\tilde I_{j,l})\;\le\;
 \Oun \sum_{l\,,\,R_{l}>{\log |I|^{-1}\over 2\log K}}
R_{l}\;\mu_{0}(\Lambda_{l})\le
 \Oun \sum_{l\,,\,R_{l}>{\log |I|^{-1}\over 2\log K}}
R_{l}\;\lambda(\Lambda_{l}) \;,
$$
and the result follows from the assumption on the exponential decay of
$\lambda(R>k)$ and Lemma II.1. 

\noindent{\bf Remark}.   Lemma II.2 implies that the density $h$ of the
measure $\mu$ with respect to the Lebesgue measure belongs to some
$L^{p}$ with $p>1$. This can also be proven directly
 using estimates similar to those in
the above proof.  Some examples of maps of the interval with neutral
fixed point are known to have an invariant measure with a density in
some $L^{p}$ (see [T]) while the bound on $\lambda(R>k)$ is only
known to be polynomial and the above proof does not work in that case.

\def\calek{{\cal E}_{k}(\epsilon)}
The proof of Theorem I.1 in the next section will require that the point
$x$ is not too rapidly recurrent. We will now prove that rapidly
recurrent  points are
exceptional with respect to the measure $\mu$.  It is convenient to
define for any integer $k$, and any positive number $\epsilon$
 the set $\calek$ by
$$
\calek=\{x\,,\,|x-f^{k}(x)|<\epsilon\}\;.
$$

\proclaim{Proposition II.3}. {There exists positive constants $C$,
$\alpha$ and $\eta<1$ such that for any integer $k$ and any $\epsilon>0$
we have
$$
\mu(\calek)\le C\left(k^{2}\epsilon^{\eta}+e^{-\alpha k}\right)\;.
$$}

\proof
We will  consider the intersection with  $\calek$ of the various
 intervals  of monotonicity of $f^{k}$. 
From formula (I.1), we have to consider
the intersection of these sets with $f^{j}(\Lambda_{l})$. 
We will start by choosing a number $\zeta>0$ such that $\beta
K^{2\zeta}<1$ and assume  first that $R_{l}<\zeta k$.
If we apply $f^{R_{l}-j}$ on $f^j(\Lambda_l)$, we land in
$\Lambda$ and we have to apply $f^{k-R_{l}+j}$.
At this point it is convenient to introduce the following construction.
Let $(s_{j})$ be a sequence of integers. We 
denote by $\Lambda_{s_{1},\,s_{2},\,\cdots\,,\,s_{r}}$ the set
$$
\Lambda_{s_{1},\,s_{2},\,\cdots\,,\,s_{r}}=\Lambda_{s_{1}}\cap
f^{-R_{s_{1}}}\Lambda_{s_{2}}\cap f^{-(R_{s_{1}}+R_{s_{2}})}\Lambda_{s_{3}}
\cap\;\cdots\;\cap
f^{-(R_{s_{1}}+\,\cdots\,+R_{s_{r-1}})}\Lambda_{s_{r}}\;.
$$
In other words, this is the subset $A$ of $\Lambda_{s_{1}}$ which is mapped
by $f^{R_{s_{1}}+\cdots+R_{s_{r-1}}}$ bijectively on $\Lambda_{s_{r}}$
with
$$
f^{R_{s_{1}}+\,\cdots\,+R_{s_{p}}}(A)\subset \Lambda_{s_{p+1}}\;.
$$
for $p=1,\;\cdots\;,\;r-1$. 

For  fixed $k$, $l$ and $j$,   we now consider all the
sets $\Lambda_{s_{1},\,\cdots\,,\,s_{r}}$ with
$R_{s_{1}}+\,\cdots\,+R_{s_{r-1}}+R_{l}-j<k$ and 
$R_{s_{1}}+\,\cdots\,+R_{s_{r}}+R_{l}-j\ge k$. Together
with $\{R>k-1-R_{l}+j\}$, this gives  a partition of $\Lambda$.
 We then construct
a partition of $f^{j}(\Lambda_{l})$ by pulling back this partition by 
$f^{R_{l}-j}$. We now consider $f^{k}$ on each atom of this partition.
Let
$$
I=I_{l,j,s_{1},\cdots,s_{r}}=f^{j}(\Lambda_{l})\cap
f^{j-R_{l}}\left(\Lambda_{s_{1},\,\cdots\,,\,s_{r}}\right)\;. 
$$

By construction, $f^{k}$ is injective on the set $I$ and we have
controled distorsion by {\bf H3}. We now prove that
 the slope of $f^{k}$ is uniformly larger
than two for $k$ large enough and $R_{s_{r}}<\zeta k$.
Let $\tilde I$ be the segment contained in $\Lambda_{l}$ which is mapped
bijectively by $f^{j}$ on $I$. By contraction and distorsion, we have
$$
|\Lambda|=|f^{R_{l}+R_{s_{1}}+\,\cdots\,+R_{s_{r}}}(\tilde I)|=
\Oun\; \left|\big(f^{R_{l}+R_{s_{1}}+\,\cdots\,+R_{s_{r}}}
\big)'_{\big|\tilde I}\right|\;
|\tilde I|
$$
$$
\le \Oun \;\left|\big(f^{R_{l}+R_{s_{1}}+\,\cdots\,+
R_{s_{r}}}\big)'_{\big|\tilde I}\right|\;
\beta^{R_{l}+R_{s_{1}}+\,\cdots\,+R_{s_{r}}}
$$
$$
\le
\Oun\;\left|f^{k'}_{\big|I}\right|\;
\beta^{k}\;K^{j}\;K^{R_{l}-j+R_{s_{1}}+\,
\cdots\,+R_{s_{r}}-k}
$$
$$
\le \Oun\; \left|f^{k'}_{\big|I}\right|\;
\beta^{k}\;K^{2\zeta k}
$$
if we assume $R_{s_{r}}<\zeta k$. The result now follows for $k$ large
enough since $\beta K^{2\zeta}<1$.

From this result it follows easily that 
if ${\calek}\cap I$ is not empty, then it is a segment
denoted below  by $J$. Assume first  that $I$ has a "large" image under
$f^{k}$, namely
$$
|f^{k}(I)|\ge\delta\;,
$$
where $\delta$ is a positive number to be chosen adequately later on.

Since $J$ is a segment, it follows easily from the definition of 
${\calek}$ that $|f^{k}(J)|\le 4\epsilon$ if we assume $|f^{k'}|>2$.
Using distorsion, we get
$$
|J|/|I|\le \Oun\epsilon/\delta\;.
$$
Using again distorsion, we get
$$
|\Lambda_l\cap f^{-j}(J)|/|\Lambda_l\cap f^{-j}(I)|\le \Oun\epsilon/\delta\;.
$$
This implies since $\mu_{0}$ is equivalent to the Lebesgue measure on
$\Lambda$ 
$$
\mu_{0}(\Lambda_l\cap f^{-j}(J))
\le \Oun{\epsilon\over \delta}\mu_{0}(\Lambda_l\cap f^{-j}(I))\;.
$$
We can now sum over all the above intervals  $I$ contained in
 $f^{j}(\Lambda_{l})$ and with "large" image. Since they are disjoint
 we get a contribution bounded above by
$\Oun(\epsilon/\delta)\mu_{0}(\Lambda_{l})$.
Summing over $j$ we get a bound $\Oun(\epsilon/\delta)R_{l}
\mu_{0}(\Lambda_{l})$.
Summing over $l$ we get finally an estimate $\Oun(\epsilon/\delta)$.
This ends the estimate in the good case of segments $I$ with "large"
images. We now have to collect the estimates for all the left-over bad cases. 

First of all we have assumed $R_{l}\le\zeta k$. We have
$$
\sum_{l\,,\,R_{l}>\zeta
k}\sum_{j=0}^{R_{l}-1}\mu_{0}\left(\Lambda_{l}\cap
f^{-j}\left({\calek}\cap f^{j}(\Lambda_{l})\right)\right)
\le \sum_{l\,,\,R_{l}>\zeta k}R_{l}\mu_{0}(\Lambda_{l})
$$
and we have a bound from Lemma II.1.

We now deal with the bad cases associated to $R_{s_{r}}$. We have by
definition 
$$
{\calek}\cap f^{j}(\Lambda_{l})=
$$
$$
\bigcup_{R_{s_{1}}+\,\cdots\,+\,R_{s_{r-1}}
<k\le R_{s_{1}}+\,\cdots\,+\,R_{s_{r}}}
 \left(I_{j,l,s_{1},\,\cdots\,,\,s_{r}}\cap{\calek}\right)
$$
$$
 \bigcup\left({\calek}\cap f^{j}(\Lambda_{l})
\cap f^{-R_{l}+j}(\{R>k-R_{l}+j\})\right)\;.
$$

We first consider the last set. We have to estimate the $\mu_{0}$
measure of
$$
\Lambda_{l}\cap f^{-j}\left({\calek}\cap f^{j}(\Lambda_{l})
\cap f^{-R_{l}+j}(\{R>k-R_{l}+j\})\right)\;.
$$
This set is obviously contained in
$$
\Lambda_{l}\cap f^{-j}\left(f^{-R_{l}+j}(\{R>k-R_{l}+j\})
\cap f^{j}(\Lambda_{l})\right)\;,
$$
which is a subset of
$$
\Lambda_{l}\cap f^{-j}\left(f^{-R_{l}+j}(\{R>(1-\zeta)k\})\right)
=\Lambda_{l}\cap f^{-R_{l}}(\{R>(1-\zeta)k\})\;,
$$
if $R_{l}\le \zeta k$ (recall that $R_{l}>j$). We get 
$$
\sum_{l,\,R_{l}<\zeta\, k}\quad\sum_{0\le j<R_{l}}
\mu_{0}\left(\Lambda_{l}\cap f^{-j}\left({\calek}\cap f^{j}(\Lambda_{l})
\cap f^{-R_{l}+j}(\{R>k-R_{l}+j\})\right)\right)
$$
$$
\le \sum_{l,\,R_{l}<\zeta\, k}\quad\sum_{0\le j<R_{l}} 
\mu_{0}\left(\Lambda_{l}\cap f^{-R_{l}}(\{R>(1-\zeta)k\})\right)
$$
$$
\le k \sum_{l,\,R_{l}<\zeta\, k}
\mu_{0}\left(\Lambda_{l}\cap f^{-R_{l}}(\{R>(1-\zeta)k\})\right)
\;.
$$
By distorsion, we have
$$
{\big|\Lambda_{l}\cap f^{-R_{l}}(\{R>(1-\zeta)k\})\big|\over
|\Lambda_{l}|}\le
$$
$$
 \Oun \left|f^{R_{l}}\left(\Lambda_{l}\cap 
f^{-R_{l}}(\{R>(1-\zeta)k\})\right)\right|
\le \Oun |\{R>(1-\zeta)k\}|\;.
$$
Therefore
$$
\sum_{l,\,R_{l}<\zeta\, k}\sum_{0\le j<R_{l}}
\mu_{0}\left(\Lambda_{l}\cap f^{-j}\left({\calek}
\cap f^{-R_{l}+j}(\{R>k-R_{l}+j\})\right)\right)
\le k \lambda(R>(1-\zeta) k)\;.
$$

We now consider the case  $R_{s_{r}}>q$ for some integer $q$.
We have
$$
\bigcup_{R_{s_{1}}+\,\cdots\,+\,R_{s_{r-1}}
<k\le R_{s_{1}}+\,\cdots\,+\,R_{s_{r}}\;,\;R_{s_{r}}>q}
 \left(I_{j,l,s_{1},\,\cdots\,,\,s_{r}}
\cap{\calek}\right)
$$
$$
\subset f^{j}(\Lambda_{l})\;\bigcap\; \left(\bigcup_{m=0}^{k}
f^{-m}(\{R>q\})\right)\;.
$$
We recall that $\{R>q\}$ is a subset of $\Lambda$. Applying $f^{-j}$ and
intersecting with $\Lambda_{l}$, we have to estimate
$$
\sum_{l\,,\, R_{l}<\zeta k}\sum_{j=0}^{R_{l}-1}\sum_{m=0}^{k}
\mu_{0}\left(\Lambda_{l}\cap f^{-j-m}(\{R>q\})\right)\;.
$$
We can now use the fact that on $\Lambda$ we have $\mu\ge\mu_{0}$.
Therefore, the above quantity is bounded by
$$
\sum_{l\,,\, R_{l}<\zeta k}\sum_{j=0}^{R_{l}-1}\sum_{m=0}^{k}
\mu\left(\Lambda_{l}\cap f^{-j-m}(\{R>q\})\right)\le
\sum_{l\,,\, R_{l}<\zeta k}\sum_{j=0}^{R_{l}-1}\sum_{m=0}^{2k}
\mu\left(\Lambda_{l}\cap f^{-m}(\{R>q\})\right)
$$
$$
\le
k\sum_{l\,,\, R_{l}<\zeta k}\sum_{m=0}^{2k}
\mu\left(\Lambda_{l}\cap f^{-m}(\{R>q\})\right)\le
k\sum_{m=0}^{2k}
\mu\left(f^{-m}(\{R>q\})\right)\le3k^{2}\mu(\{R>q\})\;.
$$

In particular, we have
$$
\mu\left(\big\{R_{s_{r}}>\zeta k\big\}\right)\le 3k^{2}\mu(\{R>\zeta k\})\;.
$$

We now have to deal with the cases $|f^{k}(I)|<\delta$. Let
$$
p=R_{l}-j+R_{s_{1}}+\,\cdots\,+R_{s_{r}}-k\;.
$$
In other words, $p$ is the number of iterations needed from $f^{k}(I)$
to reach $f^{R_{s_{r}}}(\Lambda_{s_{r}})=\Lambda$, hence 
$$
\Lambda=f^{p}(f^{k}(I))\;.
$$
Therefore
$$
|\Lambda|\le K^{p}|f^{k}(I)|\le K^{p}\delta\;,
$$
which  implies
$$
p\ge\Oun\log\delta^{-1}
$$
and therefore
$$
R_{s_{r}}\ge p\ge\Oun\log\delta^{-1}\;.
$$

We now collect all the estimates  and get
$$
\mu({\calek})\le 
$$
$$
\Oun\left({\epsilon\over \delta}+\sum_{s>\zeta k/2}\mu_{0}(R>s)+k
\mu_{0}(R>\zeta k)+k^{2}\mu(R>\zeta k)
+k^{2}\mu\big(R>\Oun\log\delta^{-1}\big)
\right)
$$
This can be expressed in terms of $\lambda$ only using Lemma II.2. 

If we assume that $\lambda(R>k)$ decays exponentially fast, namely
$$
\lambda(R>k)\le \Oun e^{-\alpha' k}
$$
for some $\alpha'>0$, we get
$$
\mu({\calek})\le 
\Oun\left({\epsilon\over\delta}+ke^{-\alpha'\zeta k}+k^{2}\delta^{\gamma}\right)
$$
for some $1>\gamma>0$. The result follows by taking the minimum with
respect to $\delta$.

In the above proof, one can avoid using the invariant measure $\mu$ in
the estimate, using instead the measure $\mu_{0}$ invariant by the map
$f^{R}$. This allows to use the same method in higher dimensional
situations. The good case corresponds to "large" enough local unstable
manifolds and give a relative bound of order $\epsilon/\delta$ which can
be integrated against the transverse measure.  The bad cases are then
handled by showing that they all correspond to large values of $R$. 

We now derive several consequences of Proposition II.3.
 Let $(E_{k})$ be the sequence of sets defined by
$$
E_{k}=\left\{y\,|\, \exists \, j\;1\le j\le (\log k)^{5}\, ,\, 
|y-f^{j}(y)|\le k^{-1}\right\}\;.
$$

\proclaim{Corollary II.4}. {There exists positive constants $C'$ and
$\beta'<1$ such that for any integer $k$
$$
\mu(E_{k})\le C' \,k^{-\beta'}\;.
$$}

\proof Note first that the estimate in Proposition II.3 is not very good
for small $k$. This can be improved as follows. We observe that since
$f$ has a slope bounded in absolute value by $K$, the inequality
$$
|f^{j}(x)-x|\le \epsilon
$$
implies
$$
|f^{2j}(x)-x|\le |f^{2j}(x)-f^{j}(x)|+|f^{j}(x)-x|\le(K^{j}+1)\epsilon\;,
$$
and more generally for any $r\ge1$
$$
|f^{rj}(x)-x|\le (K^{j}+1)^{r-1}\epsilon\;.
$$
In other words
$$
{\cal E}_{j}(\epsilon)
\subset {\cal E}_{rj}((K^{j}+1)^{r-1}\epsilon)\;.
$$
This implies together with Proposition II.3 that for any $r\ge 1$ we have
$$
\mu({\cal E}_{j}(\epsilon))\le 
\Oun \left((K^{j}+1)^{(r-1)\eta}(rj)^{2}\epsilon^{\eta}+e^{-\alpha r j}
\right)\;.
$$
Taking the minimum with respect to $r$, it follows that there are two
constants $C''>0$ and $\beta''>0$ such that for any integer $j$
$$
\mu({\cal E}_{j}(\epsilon))\le C''\epsilon^{\beta''}j^{2}\;.
$$
The  result follows by choosing $\epsilon=1/k$ and summing over $j$ from
$1$ to $(\log k)^{5}$.

\noindent{\bf Remark}. 
By a similar argument and using the Borel-Cantelli Lemma, one can show
that there is number $\rho>0$ such that the
 set of $x$ for which the event $|x-f^{k}(x)|\le k^{-\rho}$ occurs for
infinitely many $k$ is of measure zero. This would be enough for the
proof in section 3 if we use the stronger form of the decay of
correlations mentioned in the introduction. 

In order to be able to use only the weaker form of the decay of
correlations, we are going to straighten the above estimate. We will not
only control the set of points which recur too fast but also the set of
points for which a neighbor recur too fast. 
 
For positive number $\psi$ and $\rho$ to be fixed below, we define a
sequence of measurable sets $(F_{k})$  by
$$
F_{k}=\left\{x\,|\, \mu\big([x-k^{-\psi},x+k^{-\psi}]\cap E_{k^{\psi}})\ge
2\,k^{-(1+\rho)\psi}\right\}\;.
$$

\proclaim{Lemma II.5}. {The exists positive numbers $\rho$ and $\psi$
such that the  set of $x$ which belong to infinitely many
$F_{k}$ is of Lebesgue measure zero (and consequently of $\mu$ measure zero).}

\proof We will first prove that for a suitable choice of $\rho$ and
$\psi$  the sequence $(\lambda(F_{k}))$ is summable. 

Let $\bigchi_{E_{n}}$ denote the characteristic function of $E_{n}$.
We have already observed that as a consequence of Lemma II.2,
the density $h$ of $\mu$ belongs to 
$L^{p}([a,b],d\lambda)$ for some $p>1$. Therefore the 
 function $h\bigchi_{E_{n}}$
belongs also to this space. Moreover using H\"olders inequality
and Corollary II.4 its  $L^{p'}$ norm with $p'=(1+p)/2$ is bounded above
by $\Oun n^{-\vartheta}$ for some
$\vartheta>0$. We now introduce the maximal function $M_{n}$ defined by
$$
M_{n}(x)=\sup_{a>0}{1\over 2a}\int_{x-a}^{x+a}h(y)\,\bigchi_{E_{n}}(y)\,dy\;.
$$
By a well known result of Hardy and Littlewood (see [St]), this function also
belongs to $L^{p'}$ and has a norm bounded above by 
$\Oun n^{-\vartheta}$. In particular
it follows from the inequality of Chebyshev that 
$$
\lambda\left(M_{n}\ge n^{-\vartheta/2}\right)\le\Oun n^{-p'\vartheta/2}\;.
$$
In other words if $\rho=\vartheta/2$ and $\psi> 4/(p'\vartheta)$ we have
(for $k$ large enough)
$$
F_{k}\subset \left\{M_{k^{\psi}}\ge k^{-\psi\vartheta/2}\right\}
$$
which implies
$$
\lambda(F_{k})\le \Oun k^{-\psi p'\vartheta/2}\le \Oun k^{-2}\;.
$$
This last quantity is summable over $k$ and 
the result follows at once from the Borel-Cantelli Lemma.

\beginsection{III. PROOF OF THEOREM I.1.}

The strategy is inspired by the technique of extreme statistics, see for
instance [G]. We briefly explain how it works. 
Assuming  $n=pq$ with $p\approx\sqrt n$ and choosing  $s\approx
(\log n)^{2}$ we 
show that for $u_{n}=v+\log n$
$$
\proba(Z_{n}<u_n)\approx\proba(Z_{q(p+s)}<u_n)
$$
We then "dig holes" of length $s$ separating intervals of size $p$. 
Using decay of correlations we get 
$$
\proba(Z_{q(p+s)}<u_n)\approx \proba(Z_{p}<u_n)^{q}\;,
$$
and also 
$$
\proba(Z_{p}<u_n)\approx 1-p\proba(X>u_n)\;.
$$

As the reader can check, all the arguments in this section which do not
involve the results of section II work also with a
fast enough polynomial decay of correlations (with suitable choices for
$p$ and $s$).

It is convenient to use as much as possible set theoretic estimates as
presented in the next lemma. In order to alleviate the notation, we will
denote by $\{A\}$ the characteristic function of the event $A$. 

\proclaim{Lemma III.1}. { For any $k>0$ we have 
$$
\sum_{j=1}^{k}\{X_{j}>u\}\ge \{Z_{k}>u\}\ge
\sum_{j=1}^{k}\{X_{j}>u\}
-\sum_{j=1}^{k}\quad
\sum_{l\neq j\,l=1}^{k}\{X_{j}>u\}\{X_{l}>u\}\;.\eqno(III.1)
$$
}

\proof
The proof of the first inequality  is trivial, namely
 if the left hand side is zero, the right 
hand side also. On the other hand the right hand side is less than or
equal to one and the left hand side is larger than or equal to one if it
is not zero. 

For the second inequality, we have
$$
\{Z_{k}>u\}\ge
\sum_{j=1}^{k}\{X_{j}>u\}\prod_{l\neq j\,l=1}^{k}\{X_{l}<u\}
$$
i.e. if only one $X_{j}>u$ then the sup is obviously larger than $u$.
Therefore
$$
\{Z_{k}>u\}\ge\sum_{j=1}^{k}\{X_{j}>u\}
-\sum_{j=1}^{k}\{X_{j}>u\}
\left(1-\prod_{l\neq j\,l=1}^{k}\{X_{l}<u\}\right)\;.
$$
On the other hand, as in the first inequality  we have
$$
1-\prod_{l\neq j\,l=1}^{k}\{X_{l}<u\}\le\sum_{l\neq j\,l=1}^{k}\{X_{l}>u\}\;,
$$
and this implies  the lower bound.

\proclaim{Proposition III.2}. {For any integers $s,\,r,\,m,\,k,\,p\ge 0$
we have
$$
0\le\proba(Z_{r}<u)-\proba(Z_{r+k}<u)\le k\proba(X>u)\;. 
$$
and 
$$
\left|\proba(Z_{m+p+s}<u)-
\proba(Z_{m}<u)+
\sum_{j=1}^{p}
\expectation(\{X>u\}\{Z_{m}<u\}\circ f^{p+s-j})\right|\le
$$
$$
2p\sum_{j=1}^{p}\proba(\{X>u\}\{X>u\}\circ f^{j})
+s\proba(X>u)\;.
$$}

\proof  We have of course
$$
0\le \proba(Z_{r}<u)-\proba(Z_{r+k}<u)\le 
\sum_{j=0}^{k-1}\left(\proba(Z_{r+j}<u)-\proba(Z_{r+j+1}<u)\right)\;.
$$
On the other hand, for any $l\ge0$
$$
\proba(Z_{l}<u)=\proba(Z_{l+1}<u)+\proba(Z_{l}<u,\;X_{l+1}>u)
\le \proba(Z_{l+1}<u)+\proba(X_{l+1}>u)
$$
and the first  result follows by stationarity.

\noindent 
We now observe that 
$$
\{Z_{m+p+s}<u\}=\{Z_{p}<u\}\;\{Z_{s}<u\}\circ f^{p}\;
\{Z_{m}<u\}\circ f^{p+s}\;.
$$
It follows easily from this identity that
$$
|\{Z_{m+p+s}<u\}-\{Z_{p}<u\}\;\{Z_{m}<u\}\circ f^{p+s}|
\le \{Z_{s}>u\}\circ f^{p}\;.
$$

\noindent 
Therefore, using Lemma III.1 we get
$$
\big|\expectation(\{Z_{m+p+s}<u\})-
\expectation(\{Z_{p}<u\}\;\{Z_{m}<u\}\circ f^{p+s})\big|\le s\proba(X>u)\;.
$$
\noindent 
Using $\{Z_{p}<u\}=1-\{Z_{p}>u\}$, Lemma III.1 
and stationarity, we get 

$$
\left|\expectation(\{Z_{p}<u\}\;\{Z_{m}<u\}\circ f^{p+s})-
\expectation(\{Z_{m}<u\})+
\sum_{j=1}^{p}
\expectation(\{X>u\}\{Z_{m}<u\}\circ f^{p+s-j})\right|\le
$$
$$
2p\sum_{j=1}^{p}\proba(\{X>u\}\{X>u\}\circ f^{j})\;,
$$
and the result follows.

The decay of correlations is always used below in the same form, and we
present this estimate independently. It is formulated in terms
of the rate of decay $\alpha_{\omega}$ for H\"older continuous functions
of exponent $\omega$. 

\proclaim{Lemma III.3}. {For any positive number $\eta$, for any integer
$s$ and for any interval $I$ and any set $A$, we have
$$
\left|\proba(I\cap f^{-s}(A))-\proba(I)\proba(A)\right|
\le |I|^{-\omega(1+\eta)}\alpha_{\omega}(s)+\Oun |I|^{\theta(1+\eta)}\;,
$$
where $\theta$ is the number given in Lemma II.2.}

\proof The decay of correlations is formulated for H\"older continuous
functions in [Y2], and does not apply as such to characteristic
 functions. However,   if $I$ is an interval, for any number $\eta>0$
 we can find a function
$\phi$ which is non negative, satisfies $\phi\le \bigchi_{I}$, is
Lipschitz with a Lipshitz constant smaller than $|I|^{-1-\eta}$ and 
such that the support of $\bigchi_{I}(1-\phi)$ is within a distance 
$|I|^{1+\eta}$ of the boundary of $I$ 
(take for example the linear interpolation). 

We now apply the decay of correlations for functions which are H\"older
continuous with exponent $\omega$ and get
$$
\left|\int \phi\;\bigchi_{A}\circ \;f^{s}\;d\mu-
\int \phi\;d\mu\int \bigchi_{A}\;d\mu\right| 
\le |I|^{-\omega(1+\eta)}\alpha_{\omega}(s)\;.
$$
Using now Lemma II.2 we get
$$
\left|\proba(I\cap f^{-s}(A))-\proba(I)\proba(A)\right|
\le |I|^{-\omega(1+\eta)}\alpha_{\omega}(s)+\Oun |I|^{\theta(1+\eta)}\;.
$$

\noindent{\bf Remark.} The decay of correlations in [Y2] is not really
formulated for H\"older continuous functions but in term of estimates
using the function $s$. It is easy to show that any H\"older continuous
function $u$ of H\"older exponent $\omega$ satisfies these estimates.  


We now review and collect all the estimates.
We start by defining the set of full measure for which Theorem I.1
holds. This is the set of $x$ for which 
$$
\lim_{a\to 0}{1\over 2a}\mu([x-a,x+a])=h(x)
$$
and which
belong to only finitely many sets $F_{k}$ defined in section 2.
It follows from the Lebesgue differentiation
theorem applied to $\mu=h\,d\lambda$ (see for example [St]) 
and Lemma II.5 that the above two properties hold 
 for a set of full measure. 

For a fixed $v>0$ define the sequence $(u_{n})$ by
$$
u_{n}=v+\log n\;.
$$

Let $k(x)$ be the smallest integer such that $x\notin F_{j}$ for any
$j\ge k(x)$. From now on, we will assume $n>3(1+e^{-v})k(x)^{2\psi}$
 where $\psi$ is the constant appearing in Lemma II.5. 


We define the integer $p$ by $p=[\sqrt n]$ 
where $[\,\cdot\,]$ denotes the integer part. The integers $q$ and $r$
are given by the Euclidean division of $n$ by $p$, $n=pq+r$ and $0\le
r<p$. Finally we define $s=[\log n]^2$. These choices are only made for
definiteness. These choices for the numbers $p$, $q$ and $r$ are only
convenient ones. many other choices work as well.

We now  replace $\proba(Z_{n}<u_{n})$ by
$\proba(Z_{q(p+s)}<u_{n})$ and  by Proposition III.2
this produces an  error at most
$$
|\proba(Z_{n}<u_{n}) -
\proba(Z_{q(p+s)}<u_{n})|\le
qs\proba(X>u_{n})\;.
$$ 
We now estimate recursively the numbers $\proba(Z_{l(p+s)}<u_{n})$
for $0\le l\le q$. Using Lemmata III.2 and III.3 we have for any $q\ge l\ge 1$
$$
|\proba(Z_{l(p+s)}<u_{n})-(1-p\proba(X>u_{n}))
\proba(Z_{(l-1)(p+s)}<u_{n})|\le \Gamma_{n}
$$
where
$$
\Gamma_{n}=s\proba(X>u_{n})+2p\sum_{j=1}^{p}\proba
\left(\big\{X>u_{n}\big\}\cap\big\{X>u_{n}\big\}\circ f^{j}\right)
$$
$$
+p|\{X>u_{n}\}|^{-\omega(1+\eta)}\alpha_{\omega}(s)+
p\,\Oun\,|\{X>u_{n}\}|^{\theta(1+\eta)}\;.
$$
We finally get if $p\proba(X>u_n)<2$ 
$$
\left|\proba(Z_{q(p+s)}<u_{n})-(1-p\proba(X>u_{n}))^{q}\right|
\le q\Gamma_{n}\;.
$$
From Lebesgues differentiation theorem we have 
$$
\lim_{n\to\infty}pq\proba(X>u_{n})=2e^{-v}h(x)
$$
and since $s/p$ tends to zero when $n$ tends to infinity,
$$
\lim_{n\to\infty}qs\proba(X>u_{n})=0\;.
$$
A similar argument ensures $p\proba(X>u_n)<2$ for $n$ large enough. 
In order to finish the proof of Theorem I.1, we have to show that
$$
\lim_{n\to\infty}q\Gamma_{n}=0\;.
$$
If we chose $\eta$ such that $\theta(1+\eta)>2$, the result is obvious
using the exponential decay of $\alpha_{\omega}$ except for the term
$$
qp\sum_{j=1}^{p}\proba
\left(\big\{X>u_{n}\big\}\cap\big\{X>u_{n}\big\}\circ f^{j}\right)\;.
$$
Using the decay of correlations, we have easily
$$
qp\sum_{j=s}^{p}\proba
\left(\big\{X>u_{n}\big\}\cap\big\{X>u_{n}\big\}\circ f^{j}\right)
$$
$$
\le qp^{2}\proba(X>u_{n})^{2}+
qp^{2}|\{X>u_{n}\}|^{-\omega(1+\eta)}\alpha_{\omega}(s)+
qp^{2}\,\Oun\,|\{X>u_{n}\}|^{\theta(1+\eta)}\;.
$$
With the above choice of $\eta$ and the exponential decay of
$\alpha_{\omega}$, this term tends to zero when $n$ tends to infinity. 
It remains to control the part of the above sum running from $j=1$ to
$j=s-1$. 

We now define an integer $k$ (which depends on $n$) by
$$
k=\left[\big(ne^{v}/3\big)^{1/\psi}\right]\;.
$$
Recall that $n$ is large
enough so that $x$ does not belong to $F_k$. 
We now observe from the definitions that for $j\le s$ (and 
for $n$ large enough)
$$
\big\{X>u_{n}\big\}\cap
\big\{X>u_{n}\big\}\circ f^{j}
\subset [x-k^{-\psi},x+k^{-\psi}]\cap E_{k^{\psi}}\;.
$$
Since $x\notin F_{k}$ this implies 
$$
\proba\left(\big\{X>u_{n}\big\}\cap
\big\{X>u_{n}\big\}\circ f^{j}\right)\le \Oun k^{-\psi(1+\rho)}
\le \Oun n^{-(1+\rho)}\;.
$$
We finally get a bound
$$
qp\sum_{j=1}^{s}\proba
\left(\big\{X>u_{n}\big\}\cap\big\{X>u_{n}\big\}\circ f^{j}\right)
\le \Oun{qps\over n^{1+\rho}}\;,
$$
which tends to zero when $n$ tends to infinity.
This finishes the proof of Theorem I.1.

\beginsection{IV. STATISTICS OF NEAREST RECURRENCE.} 

In this section we discuss a variant of Theorem I.1 which gives the
fluctuations for the nearest return to the starting point. 
We define a sequence of real valued random variables $(X_{j})$ by 
$$
X_{j}(x)=-\log d(x,f^{j}(x))\;.
$$
We then define the sequence of random variables $(Z_{n})$ by
$$
Z_{n}(x)=\sup_{1\le j\le n}X_{j}(x)\;,
$$
and ask if the sequence of random variables $(Z_{n}-\log n)$ converges in law.
This is indeed the case under the same hypothesis as in Theorem I.1. 

\proclaim{Theorem IV.1}. {For maps of the interval satisfying the
hypothesis {\bf H1}-{\bf H4}, and such that $\lambda(R>k)$ decays
exponentially fast we have
$$
\lim_{n\to\infty}\proba\left(Z_{n}<s+\log n\right)=
\int e^{-2e^{-s}h(x)}h(x)\;dx
$$
where $h$ is the density of the invariant measure.}

Note that here also 
the normalization is related to the dimension of the measure.
In more general cases one may also expect to obtain log-normal
fluctuations as in [C.G.S.] and [K.] instead of an exponential law.

The proof is similar to that of Theorem I.1 except that we have to use
the decay of correlations to separate the initial constraint. 
We will explain in details how this can be done, and leave to the reader
to reproduce the part of the argument which is identical to the proof of
Theorem I.1.
\def\calun{{\cal U}_{n}}

\proof For a given integer $n$, let $\calun$ be the uniform partition of
the interval $[a,b]$ by intervals of length $1/n^{1+\beta'/10}$ where 
$\beta'$ is the exponent appearing in Corollary II.4 (the last segment
being of length at most this number). We fix a positive number $v$, and
from now on we will assume that $n>(1+e^{v})^2$. If $\Delta\in\calun$, we
define two intervals $\Delta^{+}$ and $\Delta^{-}$ by
$$
\Delta^{\pm}=\left\{x \,|\, d(x,\Delta)\le n^{-1}e^{-v}\pm 
n^{-1-\beta'/10}\right\}\;.
$$ 
With this notation, we have obviously
$$
\sum_{\Delta\in\calun}\proba\left(\Delta\,,\,f^{j}(\,\cdot\,)\notin
\Delta^{-}\,,\, j=1,\,\cdots\,,\,n\right)
$$
$$
\ge\proba\left(Z_{n}<v+\log n\right)\ge
$$
$$
\sum_{\Delta\in\calun}\proba\left(\Delta\,,\,f^{j}(\,\cdot\,)\notin
\Delta^{+}\,,\, j=1,\,\cdots\,,\,n\right)\;.
$$
We define $p=[n^{\theta/2}]$ ($\theta$ as given in Lemma II.2),
$s=[(\log n)^{2}]$ and let $n=(p+s)q+r$ with $0\le r<p+s$
be the division of $n$ by $p+s$. As in the first
step of the proof of Theorem I.1, we wish to replace $n$ by $q(p+s)$. 

We have obviously
$$
\left|\proba\left(\Delta\,,\,f^{j}(\,\cdot\,)\notin
\Delta^{+}\,,\, j=1,\,\cdots\,,\,n\right)-
\proba\left(\Delta\,,\,f^{j}(\,\cdot\,)\notin
\Delta^{+}\,,\, j=1,\,\cdots\,,\,(p+s)q\right)\right|\le
$$
$$
\sum_{j=q(p+s)+1}^{j=q(p+s)+r}
\proba\left(\Delta\,,\,f^{j}(\,\cdot\,)\in
\Delta^{+}\right)\;.
$$
Using decay of correlations as in Lemma III.3, we  choose 
$\eta>3/\theta$ and the above quantity is bounded by
$$
\Oun \,r\,
\left(\mu(\Delta)\mu(\Delta^{+})+n^{2\omega(1+\eta)}\alpha_{\omega}(pq)
+n^{-3-\theta}\right)\;.
$$
Using Lemma II.2, the first term is bounded by 
$$
\Oun \,r\, \mu(\Delta) n^{-\theta}\;,
$$
and since $r\le p+s\le 2n^{\theta/2}$, we can sum this quantity over
$\Delta$ and get a bound $\Oun n^{-\theta/3}$. For the two other terms,
we use the fact that the cardinality of $\calun$ is $\Oun n^{2}$ and
$\alpha_{\omega}(pq)$ decays exponentially fast in $n$. Note that this
above bounds may not apply to the last $\Delta$ in $\calun$ which may be
of size much smaller than $n^{-1-\beta'/10}$. The reader can easily
convince himself that this segment will contribute at most $\Oun
n^{-\theta(1+\beta'/10)}$ to the final result. A similar estimate holds
for the terms involving $\Delta^{-}$ instead of $\Delta^{+}$.

\def\bpqs{B^+_{p,q,s}(\Delta)} 
\def\bpqspm{B^{\pm}_{p,q,s}(\Delta)} 

For $\Delta\in\calun$ we define a set $\bpqs$ by
$$
\bpqs=\left\{x\in\Delta \,|\, f^{j}(x)\notin \Delta^{+}\; 1\le j\le (p+s)q
\right\}\;,
$$
and similarly for $B^-_{p,q,s}(\Delta)$.
From the previous bound,  we now have to estimate 
$$
\sum_{\Delta\in\calun}\proba(\bpqspm)\;.
$$
We will now eliminate the constraint $x\in\Delta$ in the definition of 
$\bpqspm$. For a fixed $v\in\real$, we assume from now on
$n$ large enough so that 
$$
e^{-({\log n})^{1/2}}> {2e^{-v}\over n}+{2\over n^{1+\beta'/10}}\;.
$$
Let
$$
G_{n}=\left\{x\;\big|\; \forall\;1\le j\le (\log n)^{2}\,,\; 
|x-f^{j}(x)|\ge e^{-\sqrt{\log n}}\right\}\;,
$$
this definition implies that if $x\in \Delta\cap G_{n}$, we have
$$
f^{j}(x)\notin \Delta^{+} \quad \hbox{\rm for} \quad
j=1,\,\cdots\,,\,\big[(\log n)^{2}\big]\;. 
$$
\def\tbpqs{\tilde B^+_{p,q,s}(\Delta)} 
Therefore, if we define $\tbpqs\supset B^+_{p,q,s}(\Delta)$ by
$$
\tbpqs=\left\{x\in\Delta \,|\, f^{j}(x)\notin \Delta^{+}\; s\le j\le (p+s)q
\right\}\;,
$$ 
we have
$$
G_{n}\cap\bpqs=G_{n}\cap\tbpqs\;.
$$
Therefore
$$
\left|\proba(\bpqs)-\proba(\tbpqs)\right|\le
\proba\big(\Delta\cap G_{n}^{c}\big)\;.
$$
The sum over $\Delta$ of this quantity is equal to $\mu(G_{n}^{c})$. 
However
$$
G_{n}^{c}\subset E_{e^{(\log n)^{1/2}}}
$$
 which implies that $\mu(G_{n}^{c})$
tends to zero when $n$ tends to infinity by Corollary II.4. It is
therefore enough to estimate $\proba(\tbpqs)$.

\def\cpqs{C^+_{p,q,s}(\Delta)} 
We now use the decay of correlations from Lemma III.3 and the estimate 
$\mu(\Delta^{+})\le\Oun n^{-\theta}$ from Lemma II.2 to replace 
$\mu(\tbpqs)$ by $\mu(\Delta)\mu(\cpqs)$ where $\cpqs$ is defined by
$$
\cpqs=\left\{x \,|\, f^{j}(x)\notin \Delta^{+}\;,\; 0\le j\le (p+s)q
\right\}\;.
$$ 
\def\calfn{{\cal F}_n^+}
The proof then proceeds following  the proof of Theorem I.1 provided
$\mu(\Delta^{+})$  is small enough, for example we can take 
$\mu(\Delta^{+}) \le |\Delta^{+}|\log n$. This is needed in order to
estimate as in section III the first part of the remainder term
$$
\Oun\;\sum_\Delta\in\calun \mu(\Delta)\left[
qs\mu(\Delta^+)+qp^2\mu(\Delta^+)^2\right]\;.
$$
 
For the second part of the remainder term, 
instead of using the sequence of sets $(F_k)$ as in section III, 
one can define for each integer $n$ a subset $\calfn$ of  $\calun$ by
$$
\calfn=\left\{\Delta\in\calun\;|\; \mu\big(\Delta^+\cap E_{[n^{2/3}]}\big)
\ge n^{-1-\beta'/2}\right\}
$$
where $\beta'$ is the constant appearing in Corollary II.4. 
We have by Corollary II.4
$$
C'n^{-2\beta'/3}\ge \mu\big(E_{[n^{2/3}]}\big)
\ge {n^{-\beta'/10}\over 4(1+e^{-v})}
 \sum_{\Delta\in\calfn}\mu\big(\Delta^+\cap
E_{[n^{2/3}]}\big)\ge {n^{-\beta'/10}\over 4(1+e^{-v})}
n^{-1-\beta'/2}\#\calfn\;,
$$
where $\#$ denotes the cardinality. The factor $n^{-\beta'/10}/4$ comes
from the fact that the sets $\Delta^+$ are not disjoint, but if we take
every other $n^{\beta'/10}e^{-v}$ such sets, we get a disjoint collection for
$n$ large enough.  Therefore
$$
\mu\left(\bigcup_{\Delta\in\calfn\;,\;\mu(\Delta^{+}) \le |\Delta^{+}|
\log n}\Delta\right)\le 
$$
$$
\mu\left(\bigcup_{\Delta\in\calfn\;,\;\mu(\Delta^{+}) \le |\Delta^{+}|
\log n}\Delta^{+}\right)\le 
\Oun n^{-1}\,\#\calfn\,\log n
\le \Oun n^{-\beta'/20}
$$
which tends to zero when $n$ tends to infinity. 
 We finally get
$$
\proba\left(Z_{n}<s+\log n\right)\ge 
\sum_{\Delta\in\calun\,,\, \mu(\Delta^{+})\le |\Delta^{+}|\log n}
\mu(\Delta)\;e^{-n\mu(\Delta^{+})}-o(1)\;.
$$
A similar upper bound follows with $\Delta^{-}$ instead of
$\Delta^{+}$, although with  an additional term, namely
$$
\proba\left(Z_{n}<s+\log n\right)\le 
\sum_{\Delta\in\calun\,,\, \mu(\Delta^{-})\le |\Delta^{-}|\log n}
\hskip -.5cm 
\mu(\Delta)\;e^{-n\mu(\Delta^{-})}+o(1)+
\sum_{\Delta\in\calun\,,\, \mu(\Delta^{-})>|\Delta^{-}|\log n
}\hskip -.5cm \mu(\Delta)
\;.
$$

By Lebesgue's derivation theorem and dominated convergence theorem, we
deduce that
$$
\lim_{n\to\infty}\sum_{\Delta\in\calun}\mu(\Delta)\; 
e^{-n\mu(\Delta^{\pm})}=
\int e^{-2e^{-v}h(x)}h(x) dx\;.
$$
It remains to control the sum of the measure of the elements $\Delta$ of
$\calun$ such that $\mu(\Delta^{\pm})>|\Delta^{\pm}|\log n$.

By Lemma II.2, it follows that the density $h$ of $\mu$ belongs to some
$L^{\sigma}$ with $\sigma>1$. Therefore, from the maximal theorem of Hardy and
Littlewood [St.] it follows that the maximal function
$$
Mh(x)=\sup_{a>0}{1\over 2a}\int_{x-a}^{x+a}h(y) dy
$$
also belongs to $L^{\sigma}$. For $\rho>0$, let $D_{\rho}$ be the set
$$
D_{\rho}=\{x\,|\, Mh(x)>\rho\}=\big\{x\,|\,\sup_{a>0}a^{-1}
\mu([x-a,x+a]>2\rho\big\}
\;.
$$
We have by Chebychev's inequality 
$$
\lambda(D_{\rho})\le \Oun\;\rho^{-\sigma}\;.
$$
We now observe that if for  a $\Delta\in \calun$ we have
$\mu(\Delta^+)\ge |\Delta^{+}|\log n$, 
then for any $y\in\Delta$ we have (for $n$ large enough)
$$
\mu([y-(1+e^{-v})n^{-1},y+(1+e^{-v})n^{-1}])
\ge\;\mu(\Delta^{+})\;\ge\; 2 e^{-v}n^{-1}\log n
$$
$$
\ge 
 \big| [y-(1+e^{-v})n^{-1},y+(1+e^{-v})n^{-1}] \big|\; (\log n)^{1/2}\;, 
$$
namely $y\in D_{(\log n)^{1/2}}$ for any $y\in\Delta$, hence $\Delta\subset
D_{(\log n)^{1/2}}$.  Therefore
$$
\lambda\left(\bigcup_{\Delta\,,\, \mu(\Delta^+)>\log n\, 
|\Delta^+|}\Delta\right)=
\sum_{\Delta\,,\, \mu(\Delta^+)>\log n\, |\Delta^+|}\lambda(\Delta)
\le \lambda\big(D_{(\log n)^{1/2}}\big)\;,
$$
which tends to zero when $n$ tends to infinity. We now use Lemma II.2 to
conclude that
$$
\lim_{n\to\infty}\mu\left(\bigcup_{\Delta\,,\, \mu(\Delta^+)>\log n\, 
|\Delta^+|}\Delta\right)=0\;.
$$
A similar argument holds
for the case of $\Delta^-$, one can also observe that $\mu(\Delta^+)\le
|\Delta^+|\log n$ implies for $n$ large enough $\mu(\Delta^-)\le
2|\Delta^-|\log n$.
This completes the proof of Theorem IV.1.\beginsection{References.}

\item{[ABST]} H.Abarbanel, R.Brown, J.Sidorowich, and L.Tsimring, 
The Analysis of
    Observed Chaotic Data in Physical Systems. Reviews of
Mod. Phys. {\bf 65}, 1331-1392 (1993). 
\medskip
\item{[BV]} V.Baladi, M.Viana.
Strong stochastic stability and rate of mixing for unimodal maps.
Ann. Sci. Ec. Norm. Super. {\bf 29},  483-517 (1996).
\medskip
\item{[CGS]} P.Collet, A.Galves and B.Schmitt. Fluctuations of
Repetition Times for Gibbsian Sources. Preprint.   
\medskip
\item{[G]} J.Galambos. {\sl The asymptotic theory of extreme order
statistics}. Wiley Series in Probability and Mathematical Statistics, 
New York (1978).
\medskip
\item{[GS]} A.Galves, B.Schmitt.
Inequalities for hitting times in mixing dynamical systems. Random Comput.
 Dyn. {\bf 5},
337-348 (1997). 
\medskip
\item{[K]} I.Kontoyiannis. Asymptotic 
Recurrence and waiting times for stationary processes. J. Theoret. Probab.
to appear.
\medskip
\item{[KN]} G.Keller,  T.Nowicki. 
Spectral theory, zeta functions and the distribution of 
periodic points for Collet-Eckmann maps.
Commun. Math. Phys. {\bf 149}, 31-69 (1992).
\medskip
\item{[Sh]} P. C. Shields. 
    The Interactions Between Ergodic Theory and Information Theory 
 IEEE
    Trans. on Information Theory. {\bf 44}, 2079-2093 (1998). 
\medskip
\item{[St]}
E.Stein. {\sl Harmonic analysis: Real-variable methods,
 orthogonality, and oscillatory integrals}.
Princeton Mathematical Series. 43. Princeton, NJ: Princeton University
 Press, 1993.
\medskip
\item{[T]} M.Thaler.
Transformations on $[0,1]$ with infinite invariant measures. 
Isr. J. Math. {\bf 46}, 67-96 (1983).
\medskip
\item{[WZW]} A. D. Wyner, J. Ziv, and A. J. Wyner. 
    On the Role of Pattern Matching in Information Theory. IEEE
    Trans. on Information Theory. {\bf 44}, 2045-2056 (1998). 
\item{[Y1]} L.S.Young. Statistical properties of dynamical systems
with some hyperbolicity. Ann. Math. {\bf 147}, 585-650 (1998).
\medskip
\item{[Y2]} L.S.Young. Recurrence times and rates of
mixing. Isr. J. Math. to appear.
\bye